\documentclass[12pt]{article}
\usepackage{amssymb}
\usepackage{latexsym,bm}
\usepackage{graphicx}
\usepackage{amsmath}
\usepackage{mathrsfs}

\date{}
\begin{document}
\title{The minimum number of detours in a connected graph of minimum degree three}
\author{\hskip -10mm Xining Liu$^{1}$, Pu Qiao$^{2,}$\thanks{Corresponding author.}, \  Xingzhi Zhan$^{1}$}
\maketitle
\footnotetext[1]{Department of Mathematics,  Key Laboratory of MEA (Ministry of Education)
 \& Shanghai Key Laboratory of PMMP, East China Normal University, Shanghai 200241, China}
\footnotetext[2]{Department of Mathematics, East China University of Science and Technology, Shanghai 200237, China}
\footnotetext[3]{E-mail addresses: {xining\_liu@126.com (X.Liu),} {pq@ecust.edu.cn (P.Qiao),} {zhan@math.ecnu.edu.cn (X.Zhan).}}

\begin{abstract}
A longest path in a graph is called a detour. Denote by $a(k,n)$ the minimum number of detours in a connected graph with minimum degree $k$ and order $n,$
and denote by $b(k,n)$ the minimum odd number of detours in such a graph. X. Zhan has posed the problem of determining $a(k,n)$ and $b(k,n).$
It is known that $a(2,n)=4$ for $n\ge 4$ and $b(2,n)=9$ for $n\ge 9.$  In this paper we prove that $a(3,n)=36$ for $n\ge 18,$ $a(k,n)\le (k!)^2$ for $n\ge k^2+2k+3$
and $b(3,n)\le 225$ for $n\ge 11.$  We also pose several related unsolved problems.
\end{abstract}


\vskip 8mm

\section{Introduction}

We consider finite simple graphs and use terminology and notations from [3] and [9]. Following Kapoor, Kronk, and Lick [8], we call a longest path in a graph $G$ a {\it detour} of $G.$ This concise term has now been widely used (e.g. [1] and [5]). In 1966, Gallai (see [6]) asked whether all detours in a connected graph share a common vertex. The answer is no in general.
It is natural to consider the number of detours in a graph.

Let $G$ be a graph. The {\it order} of $G$ is its number of vertices. The {\it size} of $G,$ denoted by $e(G),$ is the number of edges of $G.$ Thus, if $P$ is a path, then
$e(P)$ is the length of $P.$ We denote by $V(G)$ and $E(G)$ the vertex set  and edge set of a graph $G,$ respectively. For vertices $x$ and $y,$ an {\it $(x,y)$-path} is a path with endpoints $x$ and $y.$ Accordingly an {\it $(x,y)$-detour} is a detour with endpoints $x$ and $y.$ We denote by $\delta(G)$ the minimum degree of a graph $G,$ and by $N(x)$  the neighborhood of a vertex $x.$ The degree of $x$ is denoted by ${\rm deg}(x).$ If $u,\,v$ are two vertices on a path $P,$ then $P[u,v]$ denotes the subpath of $P$ with endpoints $u$ and $v.$  We denote by $K_n$ the complete graph of order $n.$

The {\it detour order} of a connected graph $G$ is defined to be the order of a detour of $G.$
A basic fact about detours [8] is that the detour order of a connected graph $G$ of order $n$ is at least ${\min}\{2\delta(G)+1,n\}.$

Let $P$ be a path in a graph $G.$ An edge $e$ in $E(P)$ is called an {\it end-edge} of $P$ if $e$ is incident with an endpoint of $P.$ An edge $e$ in $E(G)\setminus E(P)$ is called a {\it chord} of $P$ if both endpoints of $e$ lie on $P.$ A chord $e$ of a path $P$ is called a {\it boundary chord} if one endpoint of $e$ is an endpoint of $P.$
An {\it $x$-path} is a path having $x$ as an endpoint.

{\bf Notation.} $f(G)$ denotes the number of detours in a graph $G.$

The following two problems are posed in [10].

{\bf Problem 1.} Let $k$ and $n$ be integers with $2\le k\le n-2.$ Denote by $\Gamma (k,n)$ the set of connected graphs of minimum degree $k$ and order $n.$
Define
$$
a(k,n)={\rm min}\{f(G)|\, G\in \Gamma (k,n)\}.
$$
Determine $a(k,n).$

{\bf Problem 2.} Let $k,$ $n$ and $\Gamma(k,n)$ be as in Problem 1. Define
$$
b(k,n)={\rm min}\{f(G)|\, G\in \Gamma (k,n)\,\,\,{\rm and}\,\,\,f(G)\,\,\,{\rm is}\,\,\,{\rm odd}\}.
$$
Determine $b(k,n).$

First observe that for a fixed order $n,$ $a(k,n)$ is strictly increasing in $k.$ We give a proof of this fact. Let $G$ be a connected graph of minimum degree $k+1$ and order $n$ such that $f(G)=a(k+1,n)$. Let $P=x_0x_1...x_m$ be a detour of $G$ and let $e=x_0x_i$ be a boundary chord of $P$. Then $Q=P[x_{i-1},x_0]\cup x_0x_i\cup P[x_i,x_m]$ is a detour of $G.$
Since $e$ lies in a cycle, it is not a cut-edge. Thus $H=G-e$ is a connected graph of order $n$ with $k\le \delta(H)\le k+1.$ Note that $G$ and $H$ have the same detour order since $P$ is also a detour of $H$. Since every detour of $H$ is a detour of $G$ but $Q$ is not a detour of $H$, $f(H)<f(G).$ Now the condition $f(G)=a(k+1,n)$ implies that $\delta(H)=k.$
Hence $a(k,n)\le f(H)<f(G)=a(k+1,n).$

We remark that the analogues of Problems 1 and 2 for longest cycles are not very interesting, since there are uniquely hamiltonian graphs of minimum degree $4$ ([4], [7]).

It is known [10] that $a(2,n)=4$ for $n\ge 4$ and $b(2,n)=9$ for $n\ge 9.$ In this paper we prove that $a(3,n)=36$ for $n\ge 18,$ $a(k,n)\le (k!)^2$ for $n\ge k^2+2k+3$
and $b(3,n)\le 225$ for $n\ge 11.$

In Section 2 we prove the main results, and in Section 3 we pose several related unsolved problems.

\section{Main results}

We first treat hamiltonian graphs. Let $C=v_1v_2\dots v_nv_1$ be a cycle in a graph. If $i<j$, we denote by $v_i\overrightarrow{C}v_j=v_iv_{i+1}\dots v_{j}$ and $v_i\overleftarrow{C}v_j=v_iv_{i-1}\dots v_{j}$ the two $(v_i, v_j)$-paths in $C.$

{\bf Lemma 1.} {\it Let $G$ be a hamiltonian graph of order $n$ with minimum degree at least three. If $n\ge 7$ then $f(G)\ge 36.$}

{\bf Proof.} Let $C=x_1x_2...x_nx_1$ be a Hamilton cycle of $G.$ We distinguish three cases according to the value of the order $n.$

{\bf Case 1.} $7\le n\le 10.$

Denote by $\varphi(n)$ the minimum number of detours in a hamiltonian graph of minimum degree $3$ and order $n.$ A computer search shows that $\varphi(7)=48,$ $\varphi(8)=56,$ $\varphi(9)=80$ and $\varphi(10)=88.$ Hence the conclusion of Lemma 1 holds.

{\bf Case 2.} $n\ge 12.$

First, $G$ has the following $n$ detours in $C$: $x_{i}\overrightarrow{C}x_{i-1},$ $i=1, 2, ...,n.$ Here and below, the subscript of $x_j$ is to be understood modulo $n.$

Second, for every chord $x_ix_j$ $(i<j-1)$ of $C$, $G$ has the following two detours in which $x_ix_j$ is the only chord of $C$:
$$
x_{i-1}\overleftarrow{C}x_j\cup x_jx_i\cup x_i\overrightarrow{C}x_{j-1}, \quad x_{i+1}\overrightarrow{C}x_j\cup x_jx_i\cup x_i\overleftarrow{C}x_{j+1}.
$$
Since $\delta (G)\ge 3,$ the number of chords of $C$ is at least $\lceil n/2 \rceil.$ Thus $G$ has at least $2\lceil n/2 \rceil\ge n$ detours each of which contains exactly one chord of $C.$

Third, for an edge $e=x_ix_{i+1}\in E(C),$ let $x_ix_s$ and $x_{i+1}x_t$ be two chords of $C.$ Denote $\gamma(e)=\{x_ix_s,\,x_{i+1}x_t \}.$  We will construct a detour $\theta(e)$
in which $x_ix_s$ and $x_{i+1}x_t$ are the only chords of $C$ as follows. If $x_ix_s$ and $x_{i+1}x_t$ are parallel; i.e. $x_t$ lies on the path $x_{i+3}\overrightarrow{C}x_s,$
define
$$
\theta(e)=x_{s+1}\overrightarrow{C}x_i\cup x_ix_s\cup x_s\overleftarrow{C}x_t\cup x_tx_{i+1}\cup x_{i+1}\overrightarrow{C}x_{t-1}.
$$
If $x_ix_s$ and $x_{i+1}x_t$ are crossing; i.e. $x_s$ lies on the path $x_{i+2}\overrightarrow{C}x_{t-1},$
define
$$
\theta(e)=x_{t+1}\overrightarrow{C}x_i\cup x_ix_s\cup x_s\overrightarrow{C}x_t\cup x_tx_{i+1}\cup x_{i+1}\overrightarrow{C}x_{s-1}.
$$
See Figure 1.
\begin{figure}[h]
\centering
\includegraphics[width=0.5\textwidth]{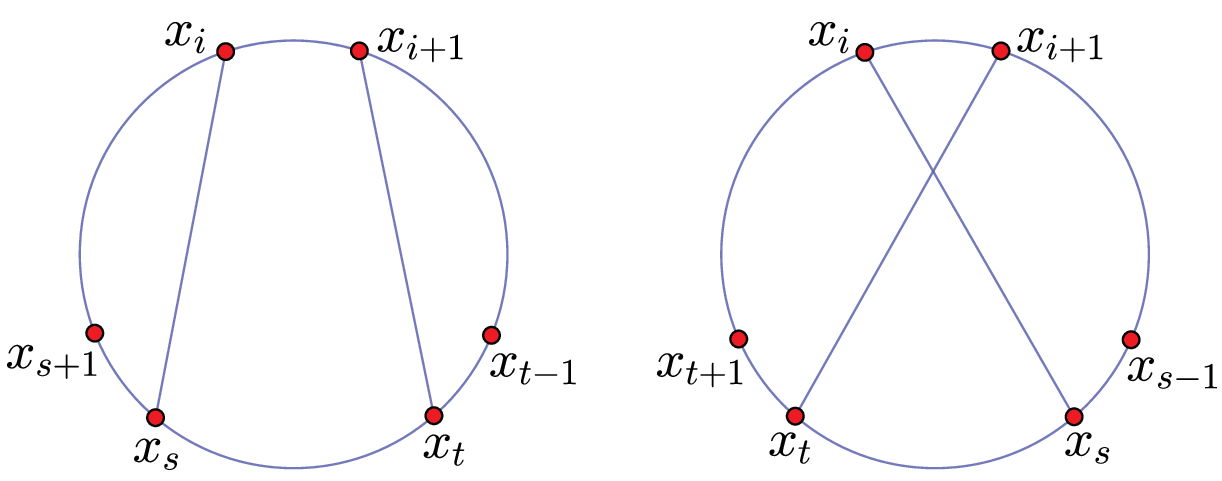}
\caption{Parallel chords and crossing chords}
\end{figure}

Observe that in the parallel case, $\theta(e)$ contains the path $x_t\overrightarrow{C}x_s$ and in the crossing case, $\theta(e)$ contains the path $x_s\overrightarrow{C}x_t.$
It is possible that for $e_1,e_2\in E(C)$ with $e_1\neq e_2,$ we have $\gamma(e_1)=\gamma(e_2).$ But in that case we still have $\theta(e_1)\neq\theta(e_2).$ Thus the $n$ edges of $C$ yield $n$ pairwise distinct detours.

Altogether, $G$ contains at least $3n\ge 36$ detours.

{\bf Case 3.} $n=11.$

As in Case 2, $G$ contains $11$ detours in $C$ and $12$ detours each of which contains exactly one chord of $C$, since $C$ has at least $\lceil 11/2 \rceil=6$ chords.
Since $\delta(G)\ge 3$ and $n=11$ is an odd number, $G$ has a vertex whose degree is at least $4.$ Without loss of generality, suppose ${\rm deg}(x_1)\ge 4.$
Then there exist two chords $x_1x_a$ and $x_1x_b$ of $C$ with $x_1$ as an endpoint. Also for every $i$ with $2\le i\le 11,$ $C$ has a chord $x_ix_{q(i)}.$  For each of the
$13$ edge-chord triples $(x_1x_2,x_1x_a,x_2x_{q(2)}),$ $(x_1x_2,x_1x_b,x_2x_{q(2)}),$ $(x_1x_{11},x_1x_a,x_{11}x_{q(11)}),$ $(x_1x_{11},x_1x_b,x_{11}x_{q(11)}),$
$(x_ix_{i+1},x_ix_{q(i)},x_{i+1}x_{q(i+1)}),i=2,\dots,10,$ $G$ contains a detour in which the two chords in that triple are the only chords of $C.$

Altogether, $G$ contains at least $11+12+13=36$ detours. \hfill$\Box$

{\bf Definition 1.} Let $P=x_1x_2...x_k$ be a path in a graph and let $x_1x_i$ be a boundary chord of $P$ with $i<k.$
From $P$ we can construct a new path of the same length that contains $x_1x_i,$ which we denote by $\Phi (P,x_1x_i)\triangleq P[x_{i-1},x_1]\cup x_1x_i\cup P[x_i,x_k].$
This classic transformation is called a {\it path exchange} [3, p.484]. We also denote the new endpoint by $\psi(P,x_1x_i)\triangleq x_{i-1}$.

We will repeatedly use the following facts.

{\bf Fact 1.}  All the neighbors of an endpoint of a detour in a graph lie on the detour.

{\bf Fact 2.}  The two endpoints of a detour in a connected nonhamiltonian graph are nonadjacent.

{\bf Fact 3.} Let $P_1$ and $P_2$ be two detours in a connected nonhamiltonian graph. Suppose that $e_1$ is a common boundary chord of $P_1$ and $P_2,$ and $e_2$ is a boundary chord of $P_1.$ Then $\Phi (P_1,e_1)\ne \Phi (P_2,e_1)$ if $P_1\ne P_2$ and $\Phi (P_1,e_1)\ne \Phi (P_1,e_2)$ if $e_1\ne e_2.$

{\bf Fact 4.} Let $P$ be an $x$-path in a graph and let $xy$ be a boundary chord of $P$ where $y$ is an internal vertex of $P.$ Then $\psi(P,xy)y\in E(P).$

We extend the concept of a chord from paths or cycles to general subgraphs.

{\bf Definition 2.} Let $H$ be a subgraph of a graph $G.$ An edge $e\in E(G)\setminus E(H)$ is called a {\it chord of $H$} if the two endpoints of $e$
lie in $H.$

{\bf Definition 3.} A {\it Wujing graph} is a graph obtained from $K_4$ by subdividing each of three edges incident with a common vertex finitely many times.

See Figure 2 for an example of a Wujing graph. Observe that every Wujing graph is traceable.

{\bf Lemma 2.} {\it Let $G$ be a connected nonhamiltonian graph of minimum degree at least three and let $H$ be a subgraph of $G$ such that $H$ is a Wujing graph.
If $H$ and $G$ have the same detour order, then $f(G)\ge 36.$}

{\bf Proof.} Let $P_i$ be an $(x,u_i)$-path in $H,$ $i=1,2,3,$ such that $V(P_i)\cap V(P_j)=\{x\}$ for $i\neq j$ and $H=P_1\cup P_2\cup P_3\cup\{u_1u_2,\, u_2u_3,\,u_3u_1\}.$

We assert that $e(P_i)\ge 2$ for each $i=1,2,3.$ Otherwise suppose $P_i=xu_i$ is an edge. Let $j\in \{1,2,3\}\setminus \{i\}.$ Then clearly $H$ contains a $(u_i, u_j)$-detour,
which is also a detour of $G.$ This contradicts Fact 2 since $G$ is connected and nonhamiltonian and $u_i$ and $u_j$ are adjacent.

Let $v_i\in N(x)$ such that $xv_i\in E(P_i)$ and let $w_i\in N(u_i)$ such that $w_iu_i\in E(P_i).$ See Figure 2.
\begin{figure}[h]
\centering
\includegraphics[width=0.3\textwidth]{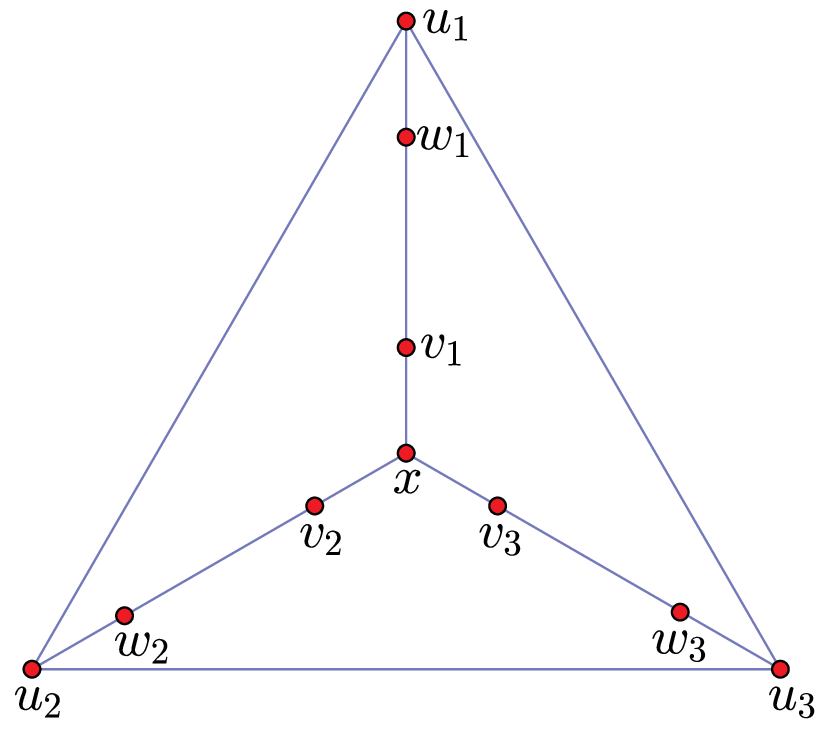}
\caption{The Wujing graph $H$}
\end{figure}

For a pair $i,j\in \{1,2,3\}$ with $i\ne j,$ let $k\in\{1,2,3\}\setminus \{i,j\}.$ There exist a $(u_i,v_j)$-detour
$$
A^{(1)}_{i,j}=P_i\cup P_k\cup u_ku_j\cup P_j[u_j,v_j]
$$
and a $(w_i,v_j)$-detour
$$
A^{(2)}_{i,j}=P_i[w_i,x]\cup P_k\cup u_ku_i\cup u_iu_j\cup P_j[u_j,v_j].
$$
Since $G$ is nonhamiltonian and for $i\ne j$ there exists a $(w_i,v_j)$-detour, by Fact 2 we deduce that $w_i$ and $v_j$ are nonadjacent. Clearly the $12$ detours $A_{i,j}^{(p)},$ $p\in \{1,2\},$ $i,j\in \{1,2,3\}$ with $i\neq j$ are pairwise distinct because of different endpoints. Since $N(w_i)\cup N(v_j)\subseteq V(A_{i,j}^{(2)})=V(H),$ ${\rm deg}(w_i) \ge 3$, and ${\rm deg}(v_j) \ge 3,$ $H$ has two chords $w_ia_i$ and  $v_jb_j,$ where $a_i\ne v_j$ and $b_j\ne w_i.$ Since $v_j$ is an endpoint of the detour $A_{i,j}^{(1)},$ $v_jb_j$ is a boundary chord of $A_{i,j}^{(1)}.$ We perform a path exchange and denote $A_{i,j}^{(3)}=\Phi (A_{i,j}^{(1)},v_jb_j).$ Since $w_i$ and $v_j$ are the two endpoints of the detour $A_{i,j}^{(2)},$ $w_ia_i$ and $v_jb_j$ are boundary chords of $A_{i,j}^{(2)}.$
We perform path exchanges and denote $A_{i,j}^{(4)}=\Phi (A_{i,j}^{(2)},w_ia_i),$ $A_{i,j}^{(5)}=\Phi (A_{i,j}^{(2)},v_jb_j)$ and $A_{i,j}^{(6)}=\Phi (A_{i,j}^{(4)},v_jb_j).$

Note that
\begin{align*}
    &E(A_{i,j}^{(1)})\setminus E(H)=E(A_{i,j}^{(2)})\setminus E(H)=\emptyset,\\
    &E(A_{i,j}^{(3)})\setminus E(H)=E(A_{i,j}^{(5)})\setminus E(H)=\{v_jb_j\},\\
    &E(A_{i,j}^{(4)})\setminus E(H)=\{w_ia_i\},\\
    &E(A_{i,j}^{(6)})\setminus E(H)=\{w_ia_i,v_jb_j\}.
\end{align*}
For $i,j,s,t\in\{1,2,3\}$ and $p,q\in\{1,\dots, 6\}$ with $i\ne j$ and $s\ne t,$ we will use the fact that if $A_{i,j}^{(p)}=A_{s,t}^{(q)},$ then $E(A_{i,j}^{(p)})\setminus E(H)=E(A_{s,t}^{(q)})\setminus E(H).$ Since $|E(A_{i,j}^{(p)})\setminus E(H)|>|E(A_{s,t}^{(q)})\setminus E(H)|,$ we have $A_{i,j}^{(p)}\ne A_{s,t}^{(q)}$ for $p\in\{3,4,5,6\}$ and $q\in \{1,2\}.$ Since $|E(A_{i,j}^{(6)})\setminus E(H)|>|E(A_{s,t}^{(q)})\setminus E(H)|,$ we have $A_{i,j}^{(6)}\ne A_{s,t}^{(q)}$ for $q\in \{3,4,5\}.$

Assume $A_{i,j}^{(3)}=A_{s,t}^{(3)}$ for $(i,j)\ne (s,t).$ Then $v_jb_j=v_tb_t.$ Since $A_{i,j}^{(1)}\ne A_{s,t}^{(1)},$ we have $A_{i,j}^{(3)}=\Phi (A_{i,j}^{(1)},v_jb_j)\ne \Phi (A_{s,t}^{(1)},v_tb_t)=A_{s,t}^{(3)},$ contradicting our assumption.

Assume $A_{i,j}^{(4)}=A_{s,t}^{(4)}$ for $(i,j)\ne (s,t).$ Then $w_ia_i=w_sa_s.$ Since $A_{i,j}^{(2)}\ne A_{s,t}^{(2)},$ we have $A_{i,j}^{(4)}=\Phi (A_{i,j}^{(2)},w_ia_i)\ne \Phi (A_{s,t}^{(2)},w_sa_s)=A_{s,t}^{(4)},$ contradicting our assumption.

Assume $A_{i,j}^{(5)}=A_{s,t}^{(5)}$ for $(i,j)\ne (s,t).$ Then $v_jb_j=v_tb_t.$ Since $A_{i,j}^{(2)}\ne A_{s,t}^{(2)},$ we have $A_{i,j}^{(5)}=\Phi (A_{i,j}^{(2)},v_jb_j)\ne \Phi (A_{s,t}^{(2)},v_tb_t)=A_{s,t}^{(5)},$ contradicting our assumption.

Assume $A_{i,j}^{(3)}=A_{s,t}^{(4)}.$ Then $v_jb_j=w_sa_s.$ Since $A_{i,j}^{(1)}\ne A_{s,t}^{(2)},$ we have $A_{i,j}^{(3)}=\Phi (A_{i,j}^{(1)},v_jb_j)\ne \Phi (A_{s,t}^{(2)},w_sa_s)=A_{s,t}^{(4)},$ contradicting our assumption.

Assume $A_{i,j}^{(3)}=A_{s,t}^{(5)}.$ Then $v_jb_j=v_tb_t.$ Since $A_{i,j}^{(1)}\ne A_{s,t}^{(2)},$ we have $A_{i,j}^{(3)}=\Phi (A_{i,j}^{(1)},v_jb_j)\ne \Phi (A_{s,t}^{(2)},v_tb_t)=A_{s,t}^{(5)},$ contradicting our assumption.

Assume $A_{i,j}^{(4)}=A_{s,t}^{(5)}.$ Then $w_ia_i=v_tb_t.$ If $i\ne t,$ then $w_i\ne v_t$ and $w_i$ is nonadjacent to $v_t.$ Thus the condition $w_ia_i=v_tb_t$ implies $i=t.$ Since $i\ne j$ and $s\ne t,$ we have $(i,j)\ne (s,t).$ Since $A_{i,j}^{(2)}\ne A_{s,t}^{(2)},$ we have $A_{i,j}^{(4)}=\Phi (A_{i,j}^{(2)},w_ia_i)\ne \Phi (A_{s,t}^{(2)},v_tb_t)=A_{s,t}^{(5)},$ contradicting our assumption.

Assume $A_{i,j}^{(6)}=A_{s,t}^{(6)}$ for $(i,j)\ne (s,t).$ Then $\{w_ia_i,v_jb_j\}=\{w_sa_s,v_tb_t\}.$ Since $A_{i,j}^{(2)}\ne A_{s,t}^{(2)},$ if $w_ia_i=w_sa_s$ and $v_jb_j=v_tb_t,$ we have 
$$
A_{i,j}^{(6)}=\Phi(\Phi(A_{i,j}^{(2)}, w_ia_i), v_jb_j)\ne \Phi(\Phi(A_{s,t}^{(2)}, w_sa_s), v_tb_t)=A_{s,t}^{(6)}.
$$
Hence $w_ia_i=v_tb_t$ and $w_sa_s=v_jb_j.$ If $i\ne t,$ then $w_i\ne v_t$ and $w_i$ is nonadjacent to $v_t.$ Thus the condition $w_ia_i=v_tb_t$ implies $i=t\ne s.$ Then we have
$$
E(A_{i,j}^{(6)})\cap \{v_ix,v_sx\}=\{v_ix\}\ne\{v_sx\}=E(A_{s,t}^{(6)})\cap \{v_ix,v_sx\},
$$
a contradiction.

Now we conclude that the $36$ detours $A_{i,j}^{(p)},$ $p\in \{1,2,\dots,6\},$ $i,j\in \{1,2,3\}$ with $i\neq j$ are pairwise distinct, implying $f(G)\ge 36.$ \hfill$\Box$

Two boundary chords $x_1x_s$ and $x_kx_t$ of a path $x_1x_2...x_k$ are said to be {\it crossing} if $s>t.$

{\bf Lemma 3.} {\it Let $G$ be a connected nonhamiltonian graph of minimum degree at least three. If $G$ has a detour with two crossing boundary chords,
then $f(G)\ge 36.$}

{\bf Proof.} Let $P=x_1x_2...x_k$ be a detour in $G$ and let $x_1x_s$ and $x_kx_t$ be two crossing boundary chords of $P$ with $s>t.$ Denote  $x=x_t$, $y=x_s$, and
$$
P_1=P[x_t,x_1]\cup x_1x_s, \quad P_2=P[x_t,x_s], \quad P_3=x_tx_k\cup P[x_k,x_s].
$$
Since $G$ is nonhamiltonian, we have $s<k,$ $t>1,$ and $s>t+1.$ It follows that $e(P_i)\ge 2$ for each $i=1,2,3.$ Note that each $P_i$ is an $(x, y)$-path and they are pairwise internally disjoint. Thus $Q\triangleq P_1\cup P_2\cup P_3$ is a $\theta$-graph whose order is equal to the detour order of $G.$
Being a $\theta$-graph, $Q$ is traceable. Let $u_i\in N(x)$ such that $xu_i\in E(P_i)$ and let $v_i\in N(y)$ such that $yv_i\in E(P_i).$
See Figure 3.
\begin{figure}[h]
\centering
\includegraphics[width=0.35\textwidth]{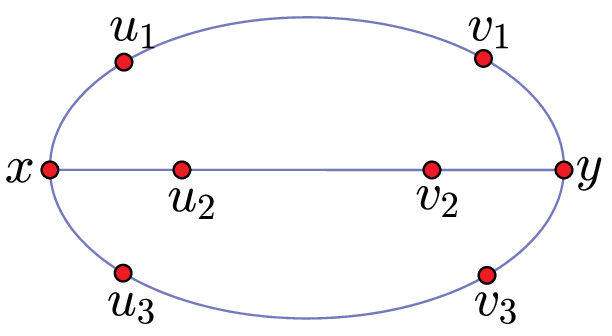}
\caption{The $\theta$-graph $Q$}
\end{figure}

Let $i, j\in \{1,2,3\}$ with $i\neq j$ and let $k\in\{1,2,3\}\setminus \{i,j\}.$
If $u_i$ and $u_j$ are adjacent, then $H_1\triangleq P_i[u_i,y]\cup P_j[u_j,y]\cup P_k \cup\{xu_i,xu_j,u_iu_j\}$ is a Wujing graph. If $v_i$ and $v_j$ are adjacent, then $H_2\triangleq P_i[v_i,x]\cup P_j[v_j,x]\cup P_k \cup\{yv_i,yv_j,v_iv_j\}$ is a Wujing graph. Since $V(H_1)=V(H_2)=V(Q)$ and $H_1,H_2,Q$ are traceable, by Lemma 2, we have $f(G)\ge 36.$

Next we assume that $u_i$ and $u_j$ are nonadjacent and $v_i$ and $v_j$ are nonadjacent for every pair $i, j\in \{1,2,3\}$ with $i\neq j.$

For a pair $i,j\in \{1,2,3\},~i\ne j$, let $k\in\{1,2,3\}\setminus \{i,j\}.$ There exists a $(u_i, v_j)$-detour
$$
A_{i,j}^{(1)}=P_i[u_i,y]\cup P_k\cup P_j[x,v_j].
$$
Clearly the six detours $A_{i,j}^{(1)},$ $i,j\in \{1,2,3\}$ with $i\neq j$ are pairwise distinct because of different endpoints. Since $N(u_i)\cup N(v_j)\subseteq V(A_{i,j}^{(1)})=V(Q),$ ${\rm deg}(u_i) \ge 3$, and ${\rm deg}(v_j) \ge 3,$ $Q$ has two chords $u_ia_i$ and $v_j b_j.$ Since $G$ is nonhamiltonian, by Fact 2 we have $a_i\ne v_j$ and $b_j\ne u_i.$ Since $u_i$ and $v_j$ are the two endpoints of the detour $A_{i,j}^{(1)},$ $u_ia_i$ and $v_jb_j$ are boundary chords of $A_{i,j}^{(1)}.$
We perform path exchanges and denote $A_{i,j}^{(2)}=\Phi (A_{i,j}^{(1)},u_ia_i),$ $A_{i,j}^{(3)}=\Phi (A_{i,j}^{(1)},v_jb_j)$ and $A_{i,j}^{(4)}=\Phi (A_{i,j}^{(2)},v_jb_j).$
Since $\tau_{i,j}^{(1)}\triangleq\psi(A_{i,j}^{(1)}, u_ia_i)$ is an endpoint of $A_{i,j}^{(2)}$ and ${\rm deg}(\tau_{i,j}^{(1)})\ge 3,$ $A_{i,j}^{(2)}$ has a boundary chord $\tau_{i,j}^{(1)} c_{i,j}.$
Since $\tau_{i,j}^{(2)}\triangleq\psi(A_{i,j}^{(1)}, v_jb_j)$ is an endpoint of $A_{i,j}^{(3)}$ and ${\rm deg}(\tau_{i,j}^{(2)})\ge 3,$ $A_{i,j}^{(3)}$ has a boundary chord $\tau_{i,j}^{(2)} d_{i,j}.$
Since $G$ is nonhamiltonian, by Fact 2 we have $c_{i,j}\ne v_j$ and $d_{i,j}\ne u_i.$ Let $u_r^+\ne x$ be the vertex such that $u_ru_r^+\in E(P_r)$ and let $v_r^-\ne y$ be the vertex such that $v_rv_r^-\in E(P_r),$ $r=1,2,3.$ Then $v_jv_j^-$ is an end-edge of $\Phi (A_{i,j}^{(2)},\tau_{i,j}^{(1)} c_{i,j})$ and $u_iu_i^+$ is an end-edge of $\Phi (A_{i,j}^{(3)},\tau_{i,j}^{(2)} d_{i,j}).$
If $a_i=y$ and $b_i=x,$ then $\tau_{i,j}^{(1)}=v_i$ and $\psi(A_{i,j}^{(2)}, v_ix)=u_k.$ If $a_j=y$ and $b_j=x,$ then $\tau_{i,j}^{(2)}=u_j$ and $\psi(A_{i,j}^{(3)}, u_jx)=v_k.$
Denote
$$
A_{i,j}^{(5)}=\begin{cases}\Phi (A_{i,j}^{(2)},\tau_{i,j}^{(1)} c_{i,j}),~~~~~~~~\,{\rm if} ~a_i\ne y~{\rm or} ~b_i\ne x,\\
\Phi(\Phi (A_{i,j}^{(2)},v_ix),u_ka_k),~{\rm if} ~a_i=y~{\rm and} ~b_i=x,\end{cases}
$$
and
$$
A_{i,j}^{(6)}=\begin{cases}\Phi (A_{i,j}^{(3)},\tau_{i,j}^{(2)} d_{i,j}),~~~~~~~~\,\,{\rm if} ~a_j\ne y~{\rm or} ~b_j\ne x,\\
\Phi(\Phi (A_{i,j}^{(3)},u_jx),v_kb_k),~~{\rm if} ~a_j=y~{\rm and} ~b_j=x.\end{cases}
$$

Note that
\begin{align*}
    &E(A_{i,j}^{(1)})\setminus E(Q)=\emptyset,\\
    &E(A_{i,j}^{(2)})\setminus E(Q)=\{u_ia_i\},\\
    &E(A_{i,j}^{(3)})\setminus E(Q)=\{v_jb_j\},\\
    &E(A_{i,j}^{(4)})\setminus E(Q)=\{u_ia_i,v_jb_j\},\\
    &E(A_{i,j}^{(5)})\setminus E(Q)=
    \begin{cases}\{u_ia_i,\tau_{i,j}^{(1)} c_{i,j}\},
    ~~~{\rm if} ~a_i\ne y~{\rm or} ~b_i\ne x,\\
    \{u_iy,v_ix,u_ka_k\},
    ~{\rm if} ~a_i=y~{\rm and} ~b_i=x,\end{cases}\\
    &E(A_{i,j}^{(6)})\setminus E(Q)=
    \begin{cases}\{v_jb_j,\tau_{i,j}^{(2)} d_{i,j}\},
    ~~~{\rm if} ~a_j\ne y~{\rm or} ~b_j\ne x,\\
    \{u_jy,v_jx,v_kb_k\},
    ~{\rm if} ~a_j=y~{\rm and} ~b_j=x.\end{cases}
\end{align*}
For $i,j,s,t\in\{1,2,3\}$ and $p,q\in\{1,\dots, 6\}$ with $i\ne j$ and $s\ne t,$ we will use the fact that if $A_{i,j}^{(p)}=A_{s,t}^{(q)},$ then $E(A_{i,j}^{(p)})\setminus E(Q)=E(A_{s,t}^{(q)})\setminus E(Q).$ Since $|E(A_{i,j}^{(p)})\setminus E(Q)|>|E(A_{s,t}^{(1)})\setminus E(Q)|,$ we have $A_{i,j}^{(p)}\ne A_{s,t}^{(1)}$ for $p\in\{2,3,4,5,6\}.$ Since $|E(A_{i,j}^{(p)})\setminus E(Q)|>|E(A_{s,t}^{(q)})\setminus E(Q)|,$ we have $A_{i,j}^{(p)}\ne A_{s,t}^{(q)}$ for $p\in\{4,5,6\}$ and $q\in \{2,3\}.$

Assume $A_{i,j}^{(2)}=A_{s,t}^{(2)}$ for $(i,j)\ne (s,t).$ Then $u_ia_i=u_sv_s.$ Since $A_{i,j}^{(1)}\ne A_{s,t}^{(1)},$ we have $A_{i,j}^{(2)}=\Phi (A_{i,j}^{(1)},u_ia_i)\ne \Phi (A_{s,t}^{(1)},u_ia_i)=A_{s,t}^{(2)},$ contradicting our assumption.

Assume $A_{i,j}^{(3)}=A_{s,t}^{(3)}$ for $(i,j)\ne (s,t).$ Then $v_jb_j=v_tb_t.$ Since $A_{i,j}^{(1)}\ne A_{s,t}^{(1)},$ we have $A_{i,j}^{(3)}=\Phi (A_{i,j}^{(1)},v_jb_j)\ne \Phi (A_{s,t}^{(1)},v_tb_t)=A_{s,t}^{(3)},$ contradicting our assumption.

Assume $A_{i,j}^{(2)}=A_{s,t}^{(3)}.$ Then $u_ia_i=v_tb_t.$ If $i\ne t,$ then $u_i\ne v_t$ and $u_i$ is nonadjacent to $v_t.$ Thus the condition $u_ia_i=v_tb_t$ implies $i=t.$ Since $i\ne j$ and $s\ne t$ we have $(i,j)\ne (s,t).$ Since $A_{i,j}^{(1)}\ne A_{s,t}^{(1)},$ we have $A_{i,j}^{(2)}=\Phi (A_{i,j}^{(1)},u_ia_i)\ne \Phi (A_{s,t}^{(1)},v_tb_t)=A_{s,t}^{(3)},$ contradicting our assumption.

Assume $A_{i,j}^{(4)}=A_{s,t}^{(4)}$ for $(i,j)\ne (s,t).$ Then $\{u_ia_i,v_jb_j\}=\{u_sa_s,v_tb_t\}.$ Since $A_{i,j}^{(1)}\ne A_{s,t}^{(1)},$ if $u_ia_i=u_sa_s$ and $v_jb_j=v_tb_t,$ we have
$$
A_{i,j}^{(4)}=\Phi(\Phi(A_{i,j}^{(1)}, u_ia_i), v_jb_j)\ne \Phi(\Phi(A_{s,t}^{(1)}, u_sa_s), v_tb_t)=A_{s,t}^{(4)}.
$$
Hence $u_ia_i=v_tb_t.$  If $i\ne t,$ then $u_i\ne v_t$ and $u_i$ is nonadjacent to $v_t.$ Thus the condition $u_ia_i=v_tb_t$ implies $i=t.$ We have
$$
E(A_{i,j}^{(4)})\cap \{u_ix,v_ty\}=\{v_ty\}\ne\{u_ix\}=E(A_{s,t}^{(4)})\cap \{u_ix,v_ty\},
$$
a contradiction.

{\bf Claim 1.} $A_{i,j}^{(5)}\neq A_{s,t}^{(5)}$ and $A_{i,j}^{(6)}\neq A_{s,t}^{(6)}$ for $(i,j)\ne (s,t).$

To the contrary assume that $A_{i,j}^{(5)}=A_{s,t}^{(5)}.$ If $\{u_iy,v_ix,u_{k_1}a_{k_1}\}=\{u_sy,v_sx,u_{k_2}a_{k_2}\}$ with $k_1\in \{1,2,3\}\setminus\{i,j\}$ and $k_2\in \{1,2,3\}\setminus\{s,t\},$ then we have $i=s.$ Thus $\{v_ix,u_{k_1}a_{k_1}\}=\{v_sx,u_{k_2}a_{k_2}\}.$ Since $(i,j)\ne (s,t)$ and $i=s,$ we have $k_1\ne k_2.$ Since $u_{k_1}\ne u_{k_2}$ and $u_{k_1}$ and $u_{k_2}$ are nonadjacent, we have $u_{k_1}a_{k_1}\ne u_{k_2}a_{k_2}.$ Since $k_2\ne s = i,$ we have $v_i\ne u_{k_2}$ and $v_i$ and $u_{k_2}$ are nonadjacent. Then $v_ix\ne u_{k_2}a_{k_2},$ a contradiction.

Now we have $\{u_ia_i,\tau_{i,j}^{(1)} c_{i,j}\}=\{u_sa_s,\tau_{s,t}^{(1)} c_{s,t}\}.$ Since $A_{i,j}^{(1)}\ne A_{s,t}^{(1)},$ if $u_ia_i=u_sa_s$ and $\tau_{i,j}^{(1)} c_{i,j}=\tau_{s,t}^{(1)} c_{s,t},$ we have
$$
A_{i,j}^{(5)}=\Phi(\Phi(A_{i,j}^{(1)}, u_ia_i), \tau_{i,j}^{(1)} c_{i,j})\ne \Phi(\Phi(A_{s,t}^{(1)}, u_sa_s), \tau_{s,t}^{(1)} c_{s,t})=A_{s,t}^{(5)}.
$$
Thus we deduce $u_ia_i=\tau_{s,t}^{(1)} c_{s,t}$ and $u_sa_s=\tau_{i,j}^{(1)} c_{i,j}.$ Note that $v_j$ is an endpoint of $A_{i,j}^{(5)}$ and $v_t$ is an endpoint of $A_{s,t}^{(5)}.$ We distinguish two cases.

{\bf Case 1.1.} $j=t.$

Then $i\ne s.$ If $a_i\in V(P_j[u_j,v_j]),$ then $\tau_{i,j}^{(1)}\in V(P_j[u_j,a_i]),$ implying $u_s\ne \tau_{i,j}^{(1)}.$ Since $u_sa_s=\tau_{i,j}^{(1)} c_{i,j},$ we have $a_s=\tau_{i,j}^{(1)},$ implying $a_s\in V(P_j[u_j,a_i]).$ Similarly, $a_i\in V(P_j[u_j,a_s]),$ a contradiction. Thus $a_i\notin V(P_j[u_j,v_j])$ and $a_s\notin V(P_j[u_j,v_j]).$
Then we have
$$
E(A_{i,j}^{(5)})\cap \{u_ix,u_sx\}=\{u_sx\}\ne\{u_ix\}=E(A_{s,t}^{(5)})\cap \{u_ix,u_sx\},
$$
a contradiction.

{\bf Case 1.2.} $j\ne t.$

Then $v_j$ and $v_t$ are the two endpoints of $A_{i,j}^{(5)}$ and $A_{s,t}^{(5)}.$ If $a_s=\tau_{i,j}^{(1)}$ and $u_s=c_{i,j},$ by Fact 4 we have $u_sv_t\in E(Q),$ a contradiction. Thus $u_s=\tau_{i,j}^{(1)},$ which implies that $s=j$ and $a_i=u_s^+=u_j^+.$ Similarly, $i=t$ and $a_s=u_i^+=u_t^+.$ Then we can verify that $A_{i,j}^{(5)}\ne A_{s,t}^{(5)},$ contradicting our assumption.

Thus we deduce $A_{i,j}^{(5)}\ne A_{s,t}^{(5)}.$ Similarly, $A_{i,j}^{(6)}\ne A_{s,t}^{(6)}.$

{\bf Claim 2.} $A_{i,j}^{(4)}\neq A_{s,t}^{(5)}$ and $A_{i,j}^{(4)}\neq A_{s,t}^{(6)}.$

To the contrary assume that $A_{i,j}^{(4)}=A_{s,t}^{(5)}.$ Then we have
$\{u_ia_i, v_jb_j\}=\{u_sa_s,\tau_{s,t}^{(1)} c_{s,t}\}.$
Since $v_tv_t^-$ is an end-edge of $A_{s,t}^{(5)}=A_{i,j}^{(4)}$ and $v_jb_j\in E(A_{i,j}^{(4)})$ we have $j\ne t.$ Then $(i,j)\ne(s,t).$ Since $A_{i,j}^{(1)}\ne A_{s,t}^{(1)},$ if $u_ia_i=u_sa_s$ and $v_jb_j=\tau_{s,t}^{(1)}c_{s,t},$ we have
$$
A_{i,j}^{(4)}=\Phi(\Phi(A_{i,j}^{(1)}, u_ia_i), v_jb_j)\ne \Phi(\Phi(A_{s,t}^{(1)}, u_sa_s), \tau_{s,t}^{(1)} c_{s,t})=A_{s,t}^{(5)}.
$$
Now we deduce $u_ia_i=\tau_{s,t}^{(1)}c_{s,t}$ and $v_jb_j=u_sa_s.$ If $j\ne s,$ then $v_j\ne u_s$ and $v_j$ is nonadjacent to $u_s.$ Thus the condition $v_jb_j=u_sa_s$ implies $j=s.$ We distinguish two cases.

If $v_jb_j=u_sa_s=u_sv_s=u_jv_j,$ then $\tau_{s,t}^{(1)}=v_s^-=v_j^-.$ Thus $a_i=\tau_{s,t}^{(1)}=v_j^-$ since $u_ia_i=\tau_{s,t}^{(1)}c_{s,t}$ and $u_i\ne v_j^-.$ Then  $v_t$ is not an endpoint of $A_{i,j}^{(4)}$  but $v_t$ is an endpoint of $A_{s,t}^{(5)},$ a contradiction.

If $v_j=u_s,$ then
$$
E(A_{i,j}^{(4)})\cap \{xv_j,yv_j\}=\{xv_j\}\ne \{yv_j\}=E(A_{s,t}^{(5)})\cap \{xv_j,yv_j\},
$$
a contradiction.

Thus we deduce $A_{i,j}^{(4)}\ne A_{s,t}^{(5)}.$ Similarly, $A_{i,j}^{(4)}\ne A_{s,t}^{(6)}.$

{\bf Claim 3.} $A_{i,j}^{(5)}\neq A_{s,t}^{(6)}.$

To the contrary assume that $A_{i,j}^{(5)}=A_{s,t}^{(6)}.$ If $\{u_iy,v_ix,u_{k_1}a_{k_1}\}=\{u_ty,v_tx,v_{k_2}b_{k_2}\}$ with $k_1\in \{1,2,3\}\setminus\{i,j\}$ and $k_2\in \{1,2,3\}\setminus\{s,t\},$ then we have $i=t$ and $u_{k_1}a_{k_1}=v_{k_2}b_{k_2}.$ If $k_1\ne k_2,$ then $u_{k_1}\ne v_{k_2}$ and $u_{k_1}$ is nonadjacent to $v_{k_2}.$  Thus the condition $u_{k_1}a_{k_1}=v_{k_2}b_{k_2}$ implies $k_1=k_2$ and $j=s.$ Then
$$
E(A_{i,j}^{(5)})\cap \{xv_j,yv_j\}=\{xv_j\}\ne \{yv_j\}=E(A_{s,t}^{(6)})\cap \{xv_j,yv_j\},
$$
a contradiction.

Thus we have $a_i\ne y$ or $b_i\ne x.$
Then $\{u_ia_i, \tau_{i,j}^{(1)} c_{i,j}\}=\{v_tb_t,\tau_{s,t}^{(2)} d_{s,t}\}.$ Note that $v_j$ is an endpoint of $A_{i,j}^{(5)}$ and $u_s$ is an endpoint of $A_{s,t}^{(6)}.$ If $v_j=u_s$, then
$$
E(A_{i,j}^{(5)})\cap \{xv_j,yv_j\}=\{xv_j\}\ne \{yv_j\}=E(A_{s,t}^{(6)})\cap \{xv_j,yv_j\},
$$
a contradiction. Thus $v_jv_j^{-}$ and $u_su_s^+$ are the two end-edges of $A_{i,j}^{(5)}=A_{s,t}^{(6)}.$ Since $u_ia_i\in E(A_{i,j}^{(5)})$ and $v_tb_t\in E(A_{s,t}^{(6)}),$
we have $i\ne s$ and $j\ne t.$ Hence $(i,j)\ne (s,t).$ Since $A_{i,j}^{(1)}\ne A_{s,t}^{(1)},$ if $u_ia_i=v_tb_t$ and $\tau_{i,j}^{(1)} c_{i,j}=\tau_{s,t}^{(2)} d_{s,t},$ we have $$
A_{i,j}^{(5)}=\Phi(\Phi(A_{i,j}^{(1)}, u_ia_i), \tau_{i,j}^{(1)} c_{i,j})\ne \Phi(\Phi(A_{s,t}^{(1)}, v_tb_t), \tau_{s,t}^{(2)} d_{s,t})=A_{s,t}^{(6)}.
$$
We deduce $u_ia_i=\tau_{s,t}^{(2)} d_{s,t}$ and $v_tb_t=\tau_{i,j}^{(1)} c_{i,j}.$ Since $v_j$ and $u_s$ are the two endpoints of $A_{i,j}^{(5)}=A_{s,t}^{(6)},$ by Fact 4 we have $u_sc_{i,j}\in E(Q)$ and $v_jd_{s,t}\in E(Q).$ Since $u_sv_t\notin E(Q)$ and $v_ju_i\notin E(Q),$ we have $v_t=\tau_{i,j}^{(1)},$ $u_i=\tau_{s,t}^{(2)},$ $b_t=c_{i,j}$ and $a_i=d_{s,t}.$ Thus $a_iv_j\in E(Q)$ and $b_tu_s\in E(Q).$ Since $a_iv_j\in E(Q)$ and $v_t=\tau_{i,j}^{(1)},$ we have $a_i=y$ and $i=t.$ Since $b_tu_s\in E(Q)$ and $u_i=\tau_{s,t}^{(2)},$ we have $b_i=b_t=x.$ This contradicts the condition that $a_i\ne y$ or $b_i\ne x.$

In summary, we have proved that $A_{i,j}^{(p)}\neq A_{s,t}^{(q)}$ for $i,j,s,t\in \{1,2,3\}$ and $p,q\in \{1,\dots,6\}$ with $i\ne j,$ $s\ne t$ and $(i,j,p)\ne (s,t,q).$ Thus we conclude that the $36$ detours $A_{i,j}^{(p)},$ $p\in \{1,2,\dots,6\},$ $i,j\in \{1,2,3\}$ with $i\neq j$ are pairwise distinct, implying $f(G)\ge 36.$\hfill$\Box$

{\bf Definition 4.} Let $P=x_0x_1...x_m$ be a detour in a graph of minimum degree at least three and let $x_0x_i$ and $x_0x_j$ be two boundary chords of $P$ with $i<j$. Let $a$
be the minimum subscript such that $x_ax_{j-1}$ is a chord of $P.$ Let $b$ be the minimum subscript such that $x_bx_{i-1}$ is a chord of $P.$ If $b\ge i+2,$ let $c$ be the
minimum subscript such that $x_cx_{b-1}$ is a chord of $P.$ If $b\le i-3,$ let $d$ be the minimum subscript such that $x_dx_{b+1}$ is a chord of $P.$ Next we define three mappings.
\begin{align*}
    &h_1(P,x_0x_i,x_0x_j)=a,\\
    &h_2(P,x_0x_i,x_0x_j)=b,\\
    &h_3(P,x_0x_i,x_0x_j)=\begin{cases}i-1,\,\,\, {\rm if}\,\,\, b=i+1,\\c,\,\,\, {\rm if}\,\,\, b\ge i+2,\\d,\,\,\, {\rm if}\,\,\, b\le i-3.\end{cases}
\end{align*}

{\bf Lemma 4.} {\it If $v$ is an endpoint of a detour in a connected nonhamiltonian graph of minimum degree at least three, then $v$ is an endpoint of six detours.}

{\bf Proof.} Let $G$ be a graph with $\delta(G)\ge 3$ and let $P=x_0x_1...x_m$ be a detour of $G$ where $x_m=v.$
 Since $N(x_0)\subseteq V(P)$ and ${\rm deg}(x_0) \ge 3,$ $P$ has two chords $x_0x_i$ and $x_0x_j$ with $2\le i<j<m.$ We will find six detours $P_1,\dots,P_6$ that have a common endpoint $x_m$.

 Recall that the functions $\Phi$ and $\psi$ are defined in Definition 1.

First let $P_1 = P,$ $P_2=\Phi (P_1,x_0x_j)$ and $P_3=\Phi (P_1,x_0x_i).$ By definition $\psi(P_1,x_0x_j)=x_{j-1}$ and $\psi(P_1,x_0x_i)=x_{i-1}.$
Since $N(x_{j-1})\subseteq V(P_2)=V(P),$ $N(x_{i-1})\subseteq V(P_3)=V(P),$ ${\rm deg}(x_{j-1}) \ge 3$ and ${\rm deg}(x_{i-1}) \ge 3,$ there exist $a=h_1(P,x_0x_i,x_0x_j)$ and $b=h_2(P,x_0x_i,x_0x_j).$ Then define $P_4=\Phi (P_2,x_{j-1}x_a)$ and $P_5=\Phi (P_3,x_{i-1}x_b).$ We distinguish two cases to
define $P_6.$ Note that since $x_bx_{i-1}$ is a chord of $P,$ we have $b\notin \{i-2, i-1, i\}.$

{\bf Case 1.} $b\ge i+1$.

By definition $\psi(P_3,x_{i-1}x_b)=x_{b-1}.$ Since $N(x_{b-1})\subseteq V(P_5)=V(P)$ and ${\rm deg}(x_{b-1}) \ge 3,$ there exists $c=h_3(P,x_0x_i,x_0x_j).$ We define $P_6=\Phi (P_5,x_{b-1}x_c).$

{\bf Case 2.} $b\le i-3$.

By definition $\psi(P_3,x_{i-1}x_b)=x_{b+1}.$ Since $N(x_{b+1})\subseteq V(P_5)=V(P)$ and ${\rm deg}(x_{b+1}) \ge 3,$ there exists $c=h_3(P,x_0x_i,x_0x_j).$ We define $P_6=\Phi (P_5,x_{b+1}x_c).$

We have found six detours that have a common endpoint $x_m$. These detours are pairwise distinct since
\begin{align*}
    &E(P_1)\setminus E(P)=\emptyset,\\
    &E(P_2)\setminus E(P)=\{x_0x_j\},\\
    &E(P_3)\setminus E(P)=\{x_0x_i\},\\
    &E(P_4)\setminus E(P)=\{x_0x_j,x_{j-1}x_a\},\\
    &E(P_5)\setminus E(P)=\{x_0x_i,x_{i-1}x_b\},\\
    &E(P_6)\setminus E(P)= \begin{cases}\{x_0x_i,x_{i-1}x_b\},~~~~~~~~~~~~{\rm if }\,\,\, b=i+1,\\\{x_0x_{b-1},x_{i-1}x_b\},\,~~~~~~~~~{\rm if }\,\,\, b\ge i+2,\, c=0,\\\{x_0x_i,x_{i-1}x_b,x_{b-1}x_c\},~~~{\rm if }\,\,\, b\ge i+2,\, c\ne0,\\\{x_{i-1}x_b, x_{b+1}x_i\},~~~~~~~~~~{\rm if }\,\,\, b\le i-3,\, c=i,\\\{x_0x_i,x_{i-1}x_b,x_{b+1}x_c\},~~~{\rm if }\,\,\, b\le i-3,\, c\ne i.\end{cases}
\end{align*}
If $b-1=i$, then $P_5\ne P_6$ since $x_ix_{i-1}\in E(P_6)\setminus E(P_5)$. This completes the proof.\hfill$\Box$

{\bf Remark.} The proof of Lemma 4 shows that each of the six detours $P_1,\dots,P_6$ contains the subpath $P[x_k, x_m]$ of $P$ where $k={\rm max}\{j,a,b,c\}.$

{\bf Lemma 5.} {\it Suppose $P=x_0x_1...x_m$ is a detour of a graph $G$ with minimum degree at least three. Let $x_0x_i,$  $x_0x_j, $ $x_mx_s$ and $x_mx_t$ be four boundary chords of $P$ with $2\le i<j<m$ and  $0<t<s\le m-2$. Denote $a=h_1(P,x_0x_i,x_0x_j),$ $b=h_2(P,x_0x_i,x_0x_j),$ $c=h_3(P,x_0x_i,x_0x_j),$ $d=h_1(P,x_mx_s,x_mx_t),$ $e=h_2(P,x_mx_s,x_mx_t),$ and $f=h_3(P,x_mx_s,x_mx_t).$
If $G$ has no detour with crossing boundary chords, then
$$
{\rm max}\{j,a,b,c\}\le {\rm min}\{t,d,e,f\}.
$$}

{\bf Proof.} Since the detour $P$ has no crossing boundary chords, $j\le t$. Define detours 
$P_1\triangleq\Phi (P,x_0x_j),$ $P_2\triangleq\Phi (P,x_0x_i),$ and $P_3\triangleq\Phi (P_2,x_{i-1}x_b).$
Since these detours have no crossing boundary chords, we have $a\le t,$ $b\le t,$ $c\le t.$

Thus for $p\in\{j,a,b,c\},$ $p\le t.$ Similarly, for $q\in\{t,d,e,f\},$ $j\le q.$

Since the detours $\Phi (P_1,x_mx_t),$ $\Phi (P_2,x_mx_t)$ and $\Phi (P_3,x_mx_t)$ have no crossing boundary chords, we deduce that  $a\le d,$
$b\le d,$ $c\le d.$ Thus for $p\in\{a,b,c\},$ $p\le d.$ Similarly, for $q\in\{d,e,f\},$ $a\le q.$

Since the detours $\Phi (P_2,x_mx_s)$ and  $P_4\triangleq\Phi (P_3,x_mx_s)$ have no crossing boundary chords, we deduce that $b\le e$ and 
$c\le e.$ Thus for $p\in\{b,c\},$ $p\le e.$ Similarly, for $q\in\{e,f\},$ $b\le q.$

Since the detour $\Phi (P_4,x_{s+1}x_e)$ has no crossing boundary chords, $c\le f.$

Finally we obtain ${\rm max}\{j,a,b,c\}\le {\rm min}\{t,d,e,f\}.$ \hfill$\Box$

{\bf Lemma 6.} {\it Let $G$ be a connected nonhamiltonian graph of minimum degree at least three. If $G$ has no detour with crossing boundary chords, then $f(G)\ge 36$.}

{\bf Proof.} Let $P=x_0x_1...x_m$ be a detour of $G$. Since $N(x_0)\subseteq V(P)$, $N(x_m)\subseteq V(P)$, ${\rm deg}(x_0) \ge 3$ and ${\rm deg}(x_m) \ge 3$, there exist four boundary chords of $P$: $x_0x_i$, $x_0x_j$, $x_mx_s$ and $x_mx_t$ where $2\le i<j<m$ and $0<t<s\le m-2$.

By Lemma 4, $G$ contains six detours $P_k,$ $k=1,\dots, 6$ that have $x_m$ as a common endpoint. By the remark after the proof of Lemma 4 and Lemma 5, we deduce that
$x_mx_s$ and $x_mx_t$ are two boundary chords of every $P_k.$ For each $k=1,\dots,6,$ applying the arguments in the proof of Lemma 4 with $P_k$, $x_mx_s$, $x_mx_t$ in place of $P,$ $x_0x_i,$ $x_0x_j,$ repeatedly, we obtain six pairwise distinct detours of $G,$ denoted by $Q_{k,r},~r=1,\dots,6.$

Denote $A_k=E(Q_{k,1})\setminus E(P)$  and $B_r=E(Q_{1,r})\setminus E(P).$ For $(k_1, r_1)\neq (k_2,r_2),$ by the proof of 
Lemma 4, $A_{k_1}\neq A_{k_2}$ or $B_{r_1}\neq B_{r_2}$. Denote
\begin{align*}
&p={\rm max}\{j,h_1(P,x_0x_i,x_0x_j),h_2(P,x_0x_i,x_0x_j),h_3(P,x_0x_i,x_0x_j)\},\\
&q={\rm min}\{t,h_1(P,x_mx_s,x_mx_t),h_2(P,x_mx_s,x_mx_t),h_3(P,x_mx_s,x_mx_t)\}.
\end{align*}
By Lemma 5, we have $p\le q.$ Let the endpoints of $Q_{k,1}$ be $\{x_{a(k)}, x_m\}$ and let the endpoints of $Q_{1,r}$ be $\{x_0, x_{b(r)}\}.$ Note that $Q_{k,r}=Q_{k,1}[x_{a(k)},x_p]\cup P[x_p,x_q]\cup Q_{1,r}[x_q,x_{b(r)}].$ Thus we conclude that
$$A_{k_1}\cap B_{r_1}=A_{k_1}\cap B_{r_2}=A_{k_2}\cap B_{r_1}=A_{k_2}\cap B_{r_2}=\emptyset
$$
and $E(Q_{k,r})\setminus E(P)=A_k\cup B_r$. Since
$$
E(Q_{k_1,r_1})\setminus E(P)=A_{k_1}\cup B_{r_1}\neq A_{k_2}\cup B_{r_2}=E(Q_{k_2,r_2})\setminus E(P),
$$
the $36$ detours $Q_{k, r},$ $1\le k,\, r\le 6$ are pairwise distinct.\hfill$\Box$

{\bf Definition 5.} {\it Adding the complete graph $K_n$ to a vertex $v$} of a graph is the operation of identifying $v$ with a vertex of $K_n.$

The operation of adding $K_4$ to a vertex is depicted in Figure 4.
\begin{figure}[h]
\centering
\includegraphics[width=0.4\textwidth]{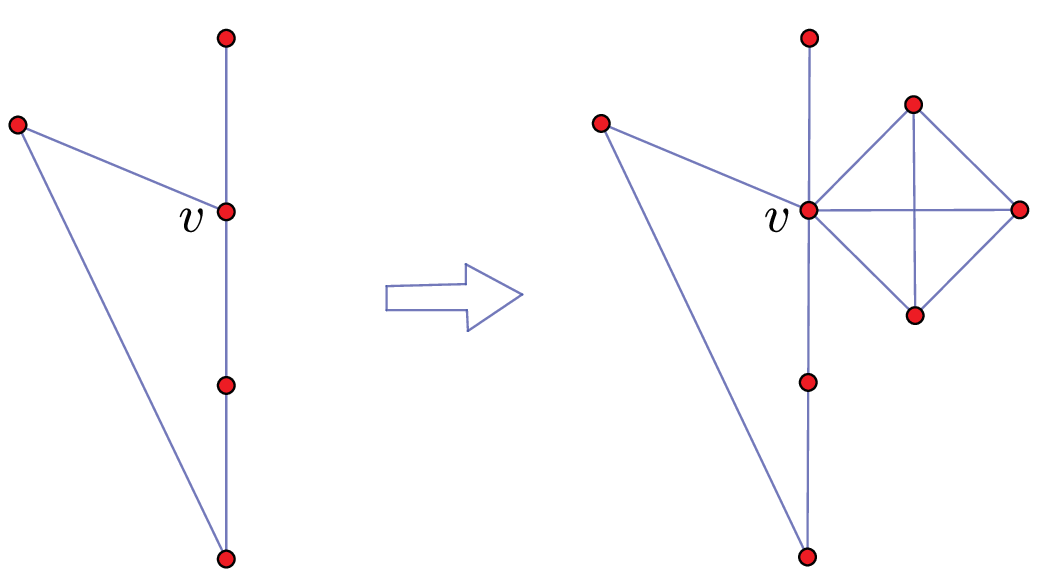}
\caption{Add $K_4$ to the vertex $v$}
\end{figure}

{\bf Definition 6.} {\it Adding the complete graph $K_n$ to an edge $e$} of a graph is the operation of first subdividing $e$ using a new vertex $v$ and then identifying $v$ with a vertex of $K_n.$

The operation of adding $K_4$ to an edge is depicted in Figure 5.
\begin{figure}[h]
\centering
\includegraphics[width=0.4\textwidth]{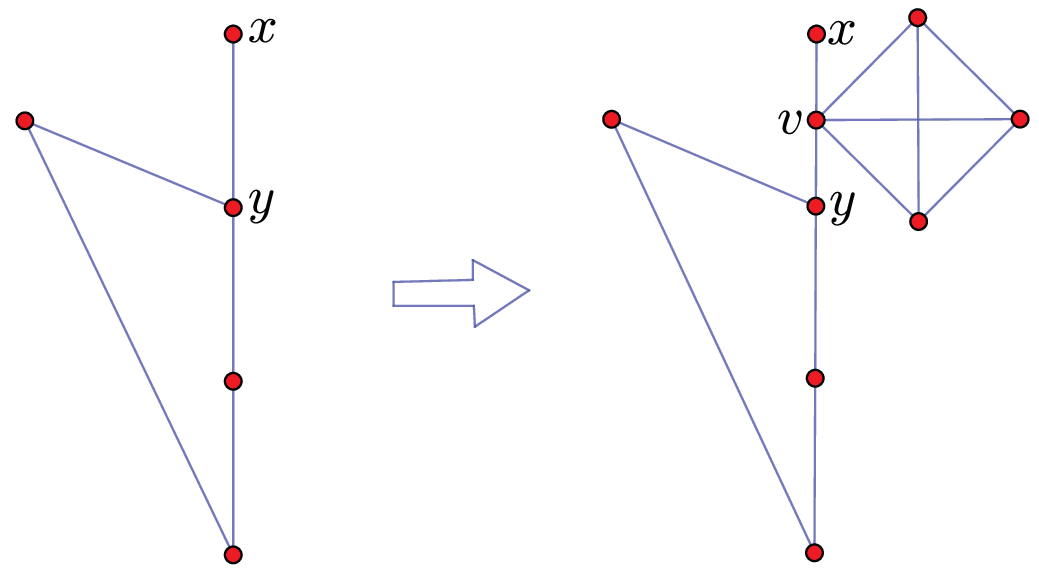}
\caption{Add $K_4$ to the edge $xy$}
\end{figure}

{\bf Theorem 7.} {\it Let $a(k,n)$ denote the minimum number of detours in a connected graph of minimum degree $k$ and order $n.$ Then
$a(3,n)=36$  if $n\ge 18$ or $n\in \{7, 8, 12, 15, 16\}.$ }

{\bf Proof.} By Lemmas 1, 3 and 6 we obtain $a(3,n)\ge 36$ for $n\ge 7.$

For every integer $n$ with  $n\in \{7,8,12,15,16\}$ or $n\ge 18,$ we construct a connected graph $G_n$ of order $n$ with $\delta(G_n)=3$ satisfying $f(G_n)=36$. 
Then the proof will be complete.

$G_7$, $G_8$, $G_{12}$, $G_{15}$ and $G_{16}$ are depicted in (a), (b), (c), (d) and (e) of Figure 6, respectively.
\begin{figure}[h]
\centering
\includegraphics[width=0.55\textwidth]{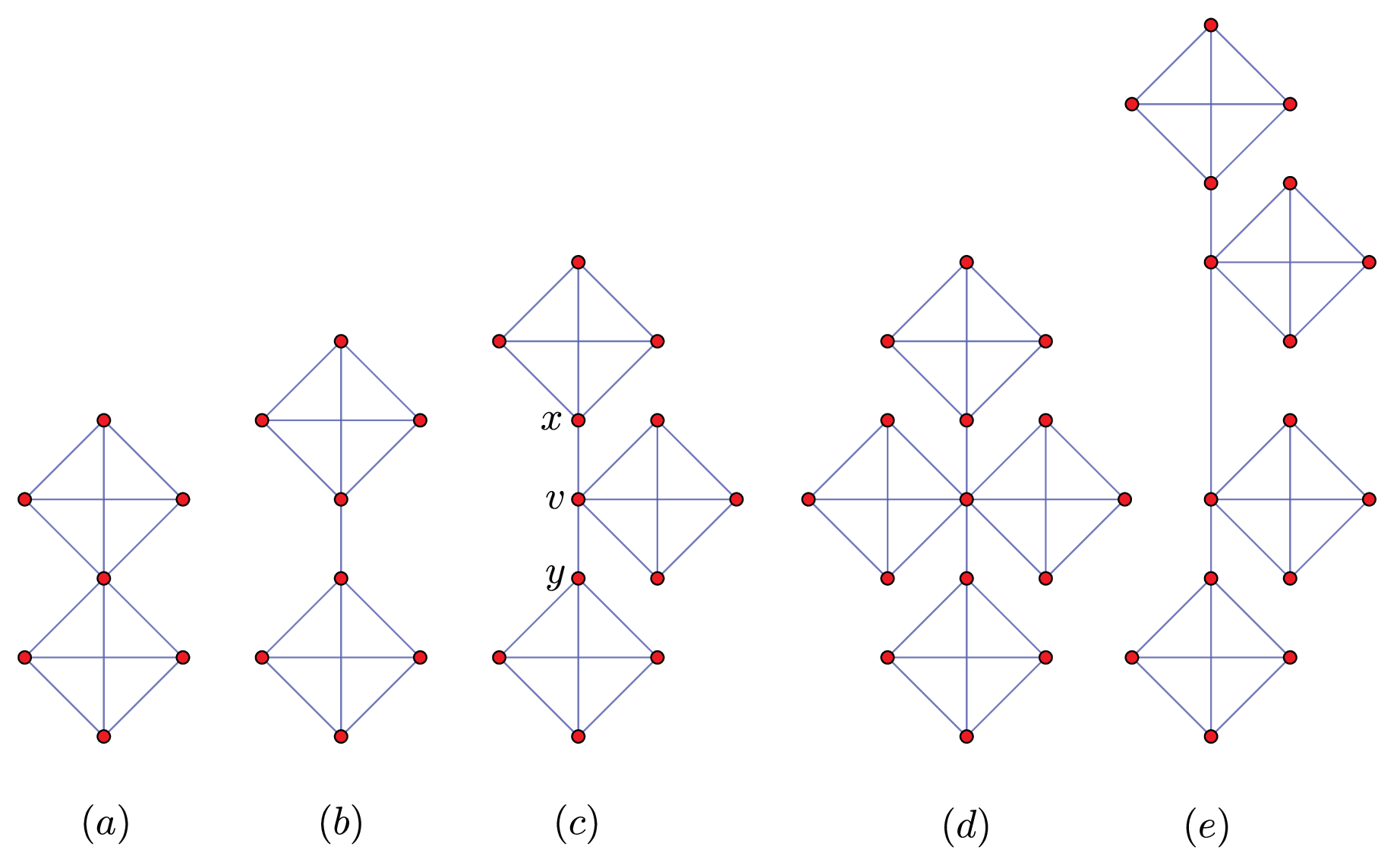}
\caption{The graphs $G_7$, $G_8$, $G_{12}$, $G_{15}$ and $G_{16}$}
\end{figure}

For every integer $n\ge 18,$ there exist nonnegative integers $p$ and $q$ with $0\le q\le 2$ satisfying $n=12+3p+4q.$ Observe that the graph $G_{12}$ has a unique central vertex which we denote by $v.$ Let the two neighbors of $v$ with degree four in $G_{12}$ be $x$ and $y.$ If $q=0,$ $G_n$ is obtained from $G_{12}$ by successively adding $K_4$ to $v$ $p$ times;
if $q=1,$ $G_n$ is obtained from $G_{12}$ by successively adding $K_4$ to $v$ $p$ times and adding $K_4$ to the edge $vx;$ if $q=2,$ $G_n$ is obtained from $G_{12}$ by successively adding $K_4$ to $v$ $p$ times, adding $K_4$ to the edge $vx$ and adding $K_4$ to the edge $vy.$ It is easy to verify that $f(G_n)=f(G_{12})=36.$ \hfill$\Box$

For example, $29 = 12+3\times3+4\times2$ and the graph $G_{29}$ is depicted in Figure 7.
\begin{figure}[h]
\centering
\includegraphics[width=0.5\textwidth]{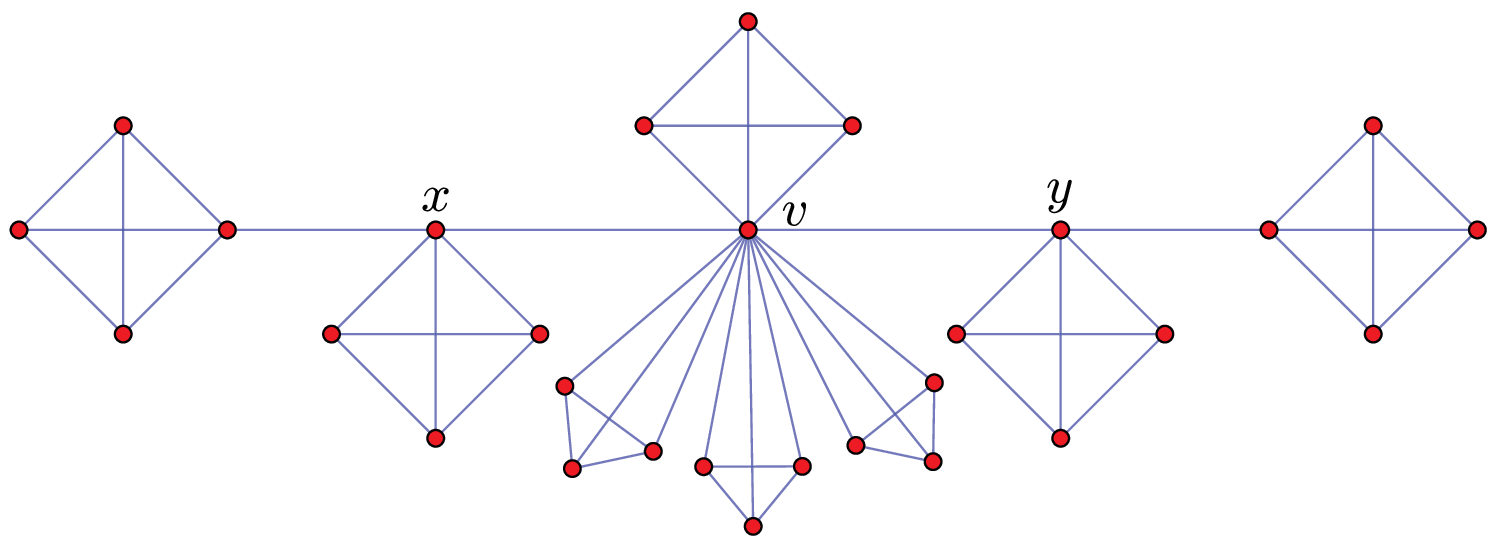}
\caption{The graph $G_{29}$}
\end{figure}

{\bf Theorem 8.} {\it Let $a(k,n)$ denote the minimum number of detours in a connected graph of minimum degree $k$ and order $n$ with $k\ge 4.$ 
If $n\in\{2k+1, 2k+2\}$ or $n=3k+3+pk+q(k+1)$ for nonnegative integers $p$ and $q,$ then $a(k,n)\le (k!)^2.$}

{\bf Proof.}  It suffices to construct a connected graph $M_n$ of order $n$ with $\delta(M_n)=k$ satisfying $f(M_n)=(k!)^2.$

For $n=2k+1,$ $M_n$ is the graph obtained by identifying a vertex of a copy of $K_{k+1}$ with a vertex of another copy of $K_{k+1}.$ For $n=2k+2,$ $M_n$ is the graph obtained by adding an edge to join a vertex of a copy of $K_{k+1}$ and a vertex of another copy of $K_{k+1}.$ For $n=3k+3,$ $M_n$ is the graph obtained by adding $K_{k+1}$ to the cut-edge
of $M_{2k+2}.$

Next suppose $n=3k+3+pk+q(k+1)$ where $p$ and $q$ are nonnegative integers and at least one of them is positive. We define $M_n$ to be the graph obtained by
first successively adding $K_{k+1}$ to the central vertex of $M_{3k+3}$ $p$ times and then successively adding $K_{k+1}$ to a cut-edge $q$ times. \hfill$\Box$

For example, $24 = 3\times 4+3+1\times 4+1\times 5$ and the graph $M_{24}$ for $k=4$ is depicted in Figure 8.
\begin{figure}[h]
\centering
\includegraphics[width=0.5\textwidth]{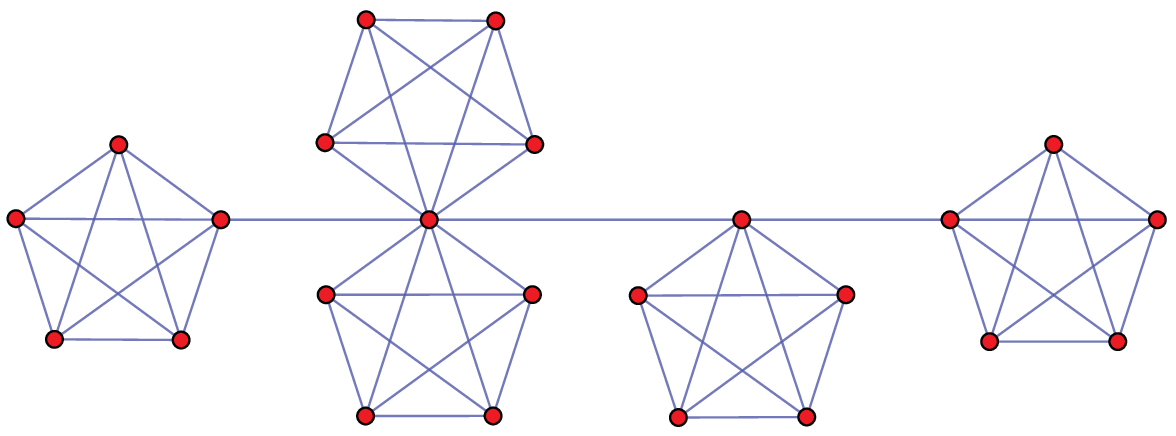}
\caption{The graph $M_{24}$ for $k=4$}
\end{figure}

{\bf Corollary 9.} {\it Let $a(k,n)$ denote the minimum number of detours in a connected graph of minimum degree $k$ and order $n$ with $k\ge 4.$
If $n\ge  k^2+2k+3,$ then $a(k,n)\le (k!)^2.$}

{\bf Proof.} Since $n\ge  k^2+2k+3=3k+3+k(k-1),$ by Theorem 8 it suffices to show that for every integer $f$ with $f\ge k(k-1),$ there exist nonnegative integers $p$ and $q$ such that $f=pk+q(k+1).$

Let $f=sk+r$ with $0\le r\le k-1.$ Then $f=(s-r)k+r(k+1)$ where $s-r\ge 0.$ \hfill$\Box$

{\bf Theorem 10.} {\it Let $b(k,n)$ denote the minimum odd number of detours in a connected graph of minimum degree $k$ and order $n.$  Then $b(3,n)\le 225$ for $n\ge 11.$}

{\bf Proof.} It suffices to construct a connected graph $R_n$ of order $n$ with $\delta(R_n)=3$ satisfying $f(R_n)=225.$

$R_{11},$ $R_{12}$ and $R_{13}$ are depicted in Figure 9.
\begin{figure}[h]
\centering
\includegraphics[width=0.55\textwidth]{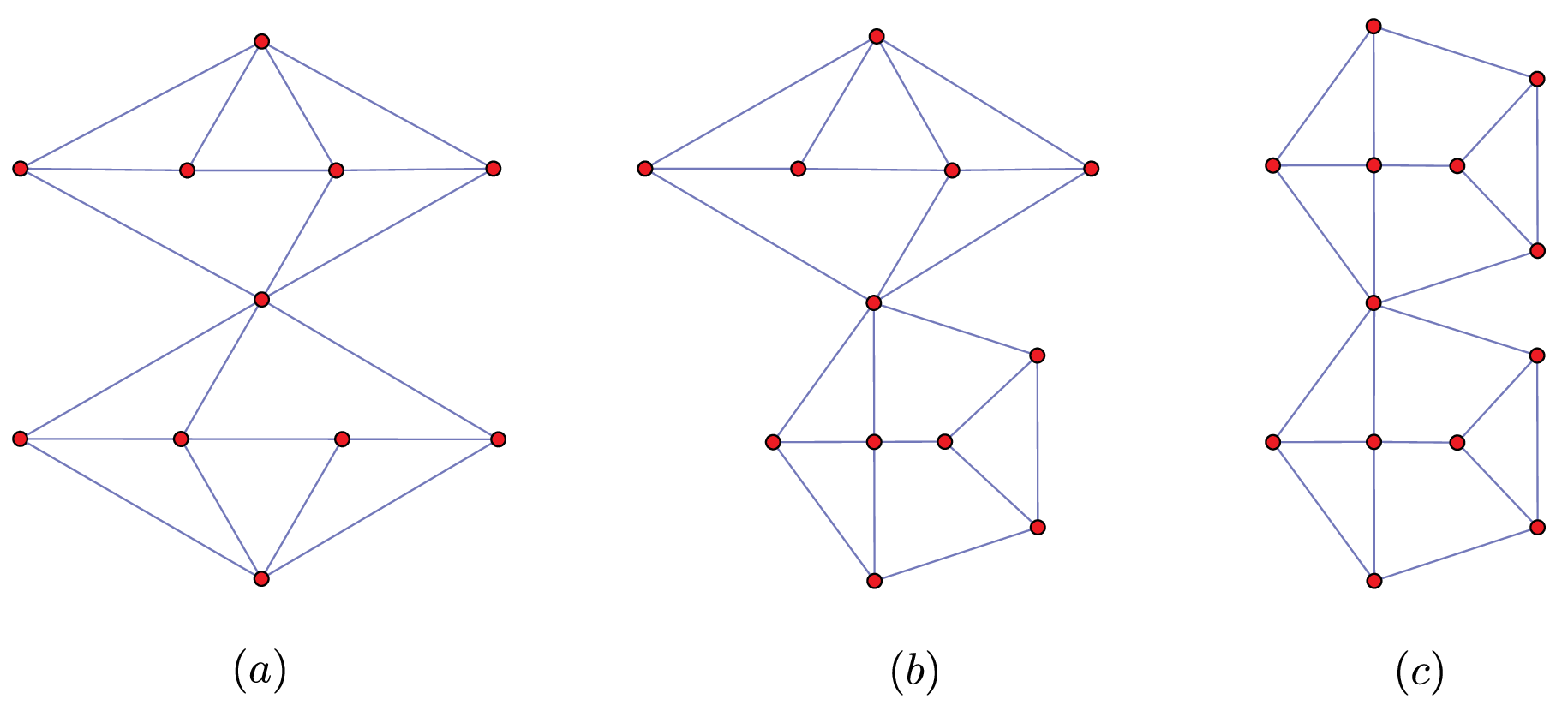}
\caption{The graphs $R_{11},$ $R_{12}$ and $R_{13}$}
\end{figure}

Next suppose $n\ge 14.$ If $n\equiv 2$ mod $3,$ then $R_n$ is obtained from $R_{11}$ by successively adding $K_4$ to the cut-vertex; if $n\equiv 0$ mod $3,$ then $R_n$ is obtained from $R_{12}$ by successively adding $K_4$ to the cut-vertex; if $n\equiv 1$ mod $3,$ then $R_n$ is obtained from $R_{13}$ by successively adding $K_4$ to the cut-vertex.

We need only to verify that each of the three graphs $R_{11},$ $R_{12}$ and $R_{13}$ has exactly $225$ detours. This can be easily done by hand or on a computer. \hfill$\Box$

\section{Unsolved problems}

Recall that $a(k,n)$ denotes the minimum number of detours in a connected graph of minimum degree $k$ and order $n,$ and $b(k,n)$ denotes the minimum odd number of detours
in such a graph; see Problems 1 and 2 in Section 1. In view of Corollary 9 and the values $a(2,n)=4$ for $n\ge 4$ and $a(3,n)=36$ for $n\ge 18,$ we pose the following conjecture.

{\bf Conjecture 3.} $a(k,n)=(k!)^2$ for $n\ge k^2+2k+3.$

The third author has posed the following problems.

{\bf Conjecture 4} (X. Zhan, February 2026). For a fixed order $n,$ $b(k,n)$ is strictly increasing in $k.$

{\bf Problem 5} (X. Zhan, February 2026). Let $k$ and $n$ be integers with $3\le k\le n-2$ and denote by $d(k,n)$ the largest integer $p$ such that if $v$ is an endpoint
of a detour in a connected graph of minimum degree $k$ and order $n,$ then $v$ is an endpoint of $p$ detours. Determine $d(k,n).$ Perhaps when the order $n$ is sufficiently large, $d(k,n)$ is independent of $n.$ It seems that $d(3,n)=8$ for $n\ge 5.$

Using the same proof ideas, we can show that the conclusion of Lemma 4 still holds when the condition ``being nonhamiltonian" is dropped. Thus $d(3,n)\ge 6.$
Next we show that $d(3,n)\le 8$ for $n\ge 5.$ It suffices to construct a connected graph $H_n$ of order $n$ with $\delta(H_n)=3$ containing a vertex $v$ that is an
endpoint of exactly eight detours of $H_n.$

$H_5,$ $H_6$ and $H_7$ are depicted in Figure 10. Let $v$ be a vertex with degree four in $H_n,$ $n=5,6,7.$ Then $v$ is an endpoint of exactly eight detours of $H_n.$
\begin{figure}[h]
\centering
\includegraphics[width=0.55\textwidth]{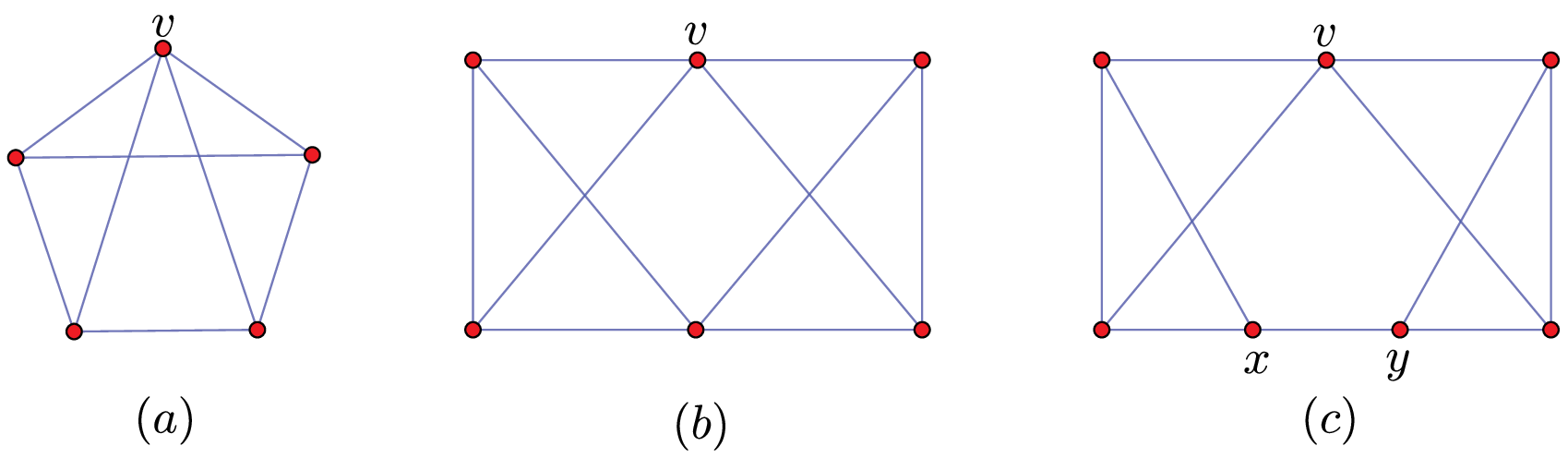}
\caption{The graphs $H_5,$ $H_6$ and $H_7$}
\end{figure}

Next suppose $n\ge 8.$ Let $x$ and $y$ be the two nonneighbors  of $v$ in $H_7.$ We define $H_n$ to be the graph obtained by subdividing the edge $xy$ with vertices $z_1,z_2,\dots,z_{n-7}$ and joining $v$ and $z_i$ for $i=1,2,\dots,n-7.$ Then it is easy to verify that $v$ is an endpoint of exactly eight detours of $H_n.$
The graph $H_{10}$ is depicted in Figure 11.
\begin{figure}[h]
\centering
\includegraphics[width=0.4\textwidth]{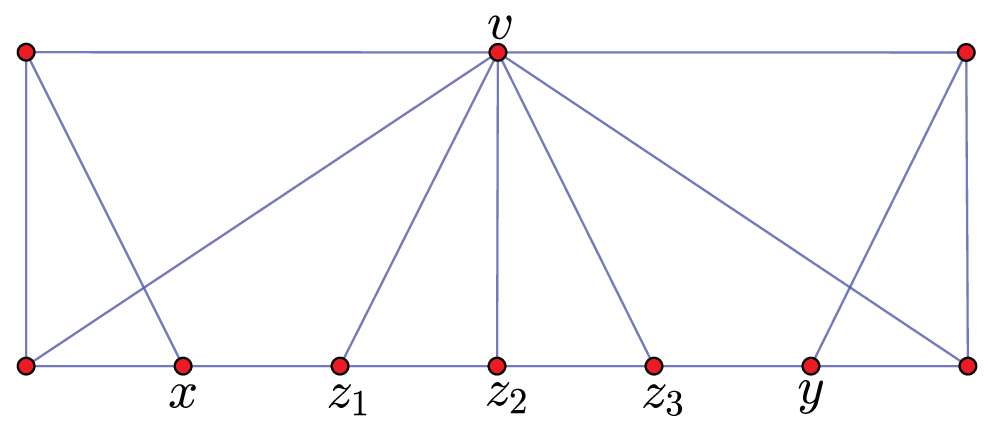}
\caption{The graph $H_{10}$}
\end{figure}

{\bf Problem 6} (X. Zhan, February 2026). Let $k$ and $n$ be integers with $3\le k\le n-2$ and denote by $q(k,n)$ the minimum number of detours in a hamiltonian graph of minimum degree $k$ and order $n.$ Determine $q(k,n).$

\vskip 5mm
{\bf Acknowledgement.} This research  was supported by the NSFC grant 12271170 and Science and Technology Commission of Shanghai Municipality
 grant 22DZ2229014.

\end{document}